\input amstex
\input amssym.tex
\documentstyle{amsppt}
\magnification=\magstep1
\hoffset=.25truein
\vsize=8.75truein

\TagsOnLeft

\define\supp{\text{supp}}
\define\Max{\text{Max}}


\define\olim{\operatornamewithlimits{\overline{lim}}}
\define\oZ{{\overline{Z}}}
\define\ep{{\epsilon}}

\define\cH{{\Cal H}}
\define\cL{{\Cal L}}
\define\cP{{\Cal P}}
\define\IC{{\Bbb C}}
\define\IN{{\Bbb N}}
\define\IP{{\Bbb P}}
\define\IR{{\Bbb R}}

\define\prf{\noindent{\bf Proof:}\quad }

\overfullrule=0pt

\topmatter
\title{Weighted Polynomials and Weighted Pluripotential Theory}\endtitle
\author{by\\  \\ 
Thomas Bloom*}\endauthor
\leftheadtext{}
\thanks{*\ Supported by an NSERC of Canada Grant.}
\endthanks

\address{Thomas Bloom, Department of Mathematics, 
University of Toronto, Toronto, Ontario  M5S 3G3 CANADA}
\endaddress
\email{bloom\@math.utoronto.ca}\endemail


\abstract
Let $E$ be a compact subset of $\IC^N$ and $w\ge 0$
a weight function on $E$ with $w>0$ on a non-pluripolar
subset of $E$. To $(E, w)$ we associate a canonical circular
set $Z\subset \IC^{N+1}$. We obtain precise relations between
the weighted pluricomplex Green function and equilibrium
measure of $(E, w)$ and the pluricomplex Green function
and equilibrium measure of $Z$. These results, combined with
an appropriate form of the Bernstein-Markov inequality,
are used to obtain asymptotic formulas for the leading
coefficients of orthonormal polynormials with respect to
certain exponentially decreasing weights in $\IR^N$.
\endabstract
\endtopmatter

\document
\baselineskip=20pt

\head  Introduction\endhead
An admissible weight on a compact set $E\subset \IC^N$
is a function $w \geq 0$ which is strictly positive on 
a non-pluripolar subset of $E$.  Associated to $(E,w)$ 
is a weighted pluripotential theory involving weighted 
polynomials, i.e, functions of the form $w^d p$ where 
$p$ is a polynomial of degrees $\le d$, a weighted pluricomplex 
Green function $V_{E,Q}$ and a weighted equilibrium 
measure $d\mu_{eq}(E,w).$  The definitions of these concepts 
are given in section 1.

In the one-dimesional case $(N=1)$ the book 
of Saff and Totik [SaTo] has many basic results. In 
the one-dimensional case, weighted polynomials arise in 
diverse problems -- approximation theory, orthognal polynomials,
random matrices, statistical physics.  For an example of recent
developments see [Dei].

In the higher dimensional case, weighted 
pluripotential theory was used
in [BL2] to obtain results on directional Tchebyshev constants
of compact sets -- the main procedure being an inductive step
from circular compact sets to a weighted problem in one less variable.

In this paper we further develop the relation between weighted
pluripotential theory on a compact set $E \subset \IC^{N}$ with
admissible (see (1.10)) weight $w$ and the
potential theory of a canonically 
associated circular
set $Z\subset \IC^{N+1}$ (defined in (2.1)).

We show that $V_Z$, the pluricomplex Green function of $Z$, and  
$d\mu_{eq}(Z)$ the equilibrium measure of $Z$, are related to 
the weighted pluricomplex Green function and the weighted 
equilibrium measure of $E$
with weight $w$.

The main results are:
\proclaim{Theorem 2.1} 
$$V_Z = (V_{E,Q}) \circ L + \log |t| \qquad {\text{for}} \ t\not=0$$
\endproclaim

\proclaim{Theorem 2.2} 
$L_*\Big({1\over 2\pi} d\mu_{eq}(Z)\Big)=d\mu_{eq}(E,w)$.
\endproclaim

Here $t$ is the first coordinate of $\IC^{N+1}$ and
$L: \big\{\IC^{N+1} - \{ (t, z)|t=0 \}\big\} \rightarrow \IC^N$ is given 
by (2.4). $L_*$ is the push-forward of measures under $L$.

Special cases of the above results may be
found in the paper of DeMarco [DeM].  In particular 
theorem 2.1 generalizes examples of section 4 of
[DeM] and theorem 2.2 generalizes lemma 2.3 of [DeM].

The advantage of considering weighted 
pluripotential theory is that (up to a limiting 
procedure described in section 5 of this paper) the
potential theory of a {\it{general}} compact circular
set in $(N + 1)$ variables may be reduced to the 
weighted case in $N$ variables. 

In section 3 we consider the Bernstein-Markov (B-M)
inequality (for the definition see (3.1)).  This 
inequality may be used to relate asymptotics of
orthonormal polynomials with respect to a measure $\mu$ on $E$
to potential theoretic invariants of $E$.  We introduce a weighted 
version of the B-M
inequality (see(3.2)).  We show (theorem 3.1) that the weighted
B-M inequality holds on $E$ with weight $w$ and measure $\mu$
if and only if the B-M 
inequality holds for an associated measure on
$Z$. 
Then we give the following cases where the B-M
inequality holds.

\proclaim{Corollary 3.1}
$d\mu_{eq}(E,w)$ for $E$ regular and weight $w$.
\endproclaim

\proclaim{Theorem 3.2} 
$(E, w, \sigma)$ where $E\subset \IR^N$ and $(E, \sigma)$
staisfies the B-M inequality.
\endproclaim

In section 4 we obtain asymptotics for the 
leading coefficients of orthonormal polynomials with 
respect to certain exponentially decreasing measures
on $\IR^n$. As in the known procedure in the one-variable
case, we first scale the problem to obtain a
problem on the asymptotics of weighted polynomials.
Using the weighted B-M inequality (a special case
of theorem 3.2) gives the asymptotics (see example 4.1
and equation (4.24)).

\head 1.\ Preliminaries\endhead
Let $E$ be a bounded subset of $\IC^N$.
The pluricomplex Green function of $E$ is defined
by
$$
V_E (z):=\sup \{u(z)|u\in \cL, u\le 0\ {\text{on}}\ E\}\ {\text{where}}
\tag 1.1 $$ 
$$
\cL=\{u|u \ \ {\text{is\ plurisubharmonic\ (PSH)\ on}}\ \IC^N,
u (z) \le \log^+ |z| +C\}.
\tag 1.2$$
is the Lelong class of PSH functions
of logarithmic growth. (We use the notation
$$
|z|:=\left(\sum^N_{i=1} |z_i|^2\right)^{1/2}\ {\text{for}}\
z=(z_1,\cdots, z_N) \in \IC^N).$$

A set $E\subset \IC^N$ is said to be
{\it{pluripolar}} if for all points $a\in E$, there
is a neighborhood $U$ of $a$ and a function
$v$ which is PSH on $U$ such that
$E\cap U\subset\{z\in U| v  (z)=-\infty\}$. A property of
a set $E$ is said to hold {\it{quasi-everywhere}} (q.e.)
if there is a pluripolar set $P\subset E$ and the
property holds at all points of $E\setminus P$.

For $G$ an open subset of $\IC^N$ and
$f$ a real-valued function on $G$, we let $f^*$
denote its uppersemicontinuous (u.s.c.) 
regularization, defined by
$$f^*(z)=\olim_{\xi\to z} f(\xi)\qquad {\text{for}}\ z\in G.$$
$V^*_E\in \cL$ if and only if $E$ is non-pluripolar [K].
For $E$ non-pluripolar, the equilibrium measure
of $E$ is defined by
$$d\mu_{eq}=d\mu_{eq} (E):=(dd^c V^*_E)^N
\tag 1.3$$
where $(dd^c)^N$ is the complex Monge-Amp\`ere operator.
$d\mu_{eq}$ is a positive Borel measure of total
mass $(2\pi)^N$ and with $\supp (d\mu_{eq})\subset \overline E$. [K].

In the case that $E$ is a compact set, it
is a result of Siciak and Zaharyuta ([K], theorem 5.1.7)
that 
$$V_E(z)=\log \phi_E (z)
\tag 1.4$$
where
$$\phi_E (z)=\sup\{|p(z)|^{1\over{\deg (p)}}| p\
{\text{is\ a\ holomorphic\ polynomial}},\ \deg(p)\ge 1\
{\text{and}}\ \|p\|_E \le 1\}.
\tag 1.5$$
It follows that $V_E(z)=V_{\hat E}(z)$ where $\hat E$ denotes
the polynomially convex hull of $E$.

We will also use the class $\cH$
of logarithmically homogeneous PSH functions
on $\IC^N$ defined by
$$\cH:=\{u\in \cL | u(tz)=u(z)+\log |t|\ {\text{for\ all}}\
z\in \IC^N, t\in \IC\}.
\tag 1.6$$
For $E$ a bounded set in $\IC^n$, we define
$$H_E(z):=\sup\{u(z)|u\in\cH,\ u\le 0\ {\text{on}} E\}.
\tag 1.7$$
For $E$ compact (see [Si2]) we have
$$H_E(z)=\log \psi_E (z)\ {\text{where}}
\tag 1.8$$
$$\psi_E(z)=\sup\{|p(z)|^{1\over \deg(p)}|p\ 
{\text{is\ a\ homogeneous\ holomorphic\ polynomial}},\ 
\deg p\ge 1\ {\text{and}} \|p\|_E \le 1\}.
\tag 1.9$$
For $E$ a compact set in $\IC^N$, an {\it{admissible
weight function}} is a real-valued function on $E$
satisfying
$$\aligned i)\quad & w\ge 0\\
	 ii)\quad & w\ {\text{is\ u.s.c.}}\\
	 iii)\quad & \{z\in E|w(z)>0\}\ {\text{is\ non-pluripolar.}}
	 \endaligned\tag 1.10
$$
In particular, if $E$ admits on admissible weight
function then $E$ itself is non-pluripolar.

There is a ``weighted" version of the
pluricomplex Green function (see [S.1], [SaTo, appendix B])
defined as follows: let
$$Q:=-\log w\tag 1.11$$
Then $Q$ is lowersemicontinuous (l.s.c.) on $E$.
The weighted pluricomplex Green function of $E$ with 
weight $w$ is defined by
$$V_{E, Q}:=\sup\{u(z)|u\in\cL,\ u\le Q\ {\text{on}}\ E\}
\tag 1.12$$

The weighted equilibrium measure of $E$ is defined by
$$d\mu_{eq} (E, w):=(dd^c V^*_{E, Q})^N.\tag 1.13$$
It is a positive Borel measure with supp $(d\mu_{eq} (E, w))\subset E$.
and total mass $(2\pi)^N$.

A weighted polynomial on $E$ is defined to
be a function of the form $w^d p$ where
$d$ is an integer $\ge 0$ and $p$ is a holomorphic
polynomial of degree $\le d$. Note that if $\| w^d p\|E \le 1$
then ${1\over d}\log |p(z)| \le Q(z)$ on $E$ and since
${1\over d}\log |p(z)| \in \cL$
we have, ${1\over d} \log |p(z)| \le V_{E, Q} (z)$ for
all $z\in \IC^N$.

It is known (see [Si1] or [SaTo], appendix B) that
$$V_{E, Q} (z)=\log \phi_{E, Q}(z)\ \ {\text{where}}\tag 1.14$$
$$\phi_{E,Q}(z)=\sup\{|p(z)|^{1\over d}| \|w^d p\|_E\le 1,
\deg p\ge 1\ {\text{and}}\
w^d p\ {\text{is\ a\ weighted\ polynomial}}\}.\tag 1.15$$

A set $E$ is defined to be {\it{regular}} if $V_E$ is
continuous on $\IC^N$. A set $E$ is defined ([Si1]) to
be {\it{locally regular at a point}} $a\in \overline E$ if for
each $r>0$, $V_{E\cap B(a, r)}$ is continuous at $a$. Here
$B(a, r):=\{z\in \IC^N\bigm|  |z-a|\le r\}$ denotes the ball
center $a$, radius $r$. It is sufficient, for $E$ to be
locally regular at $a$, that $V_{E\cap B(a, r)}$ be continuous at $a$
for all $r>0$ sufficiently small.

$E$ is said to be {\it{locally regular}} if it is
locally regular at each point of $\overline E$.

For $E$ compact and locally regular and $w$ a continuous admissible
weight function on $E$ then $V_{E, Q}$ is continuous [Si1].

For $u\in \cL$ we define its Robin function $\rho_u$
by
$$\rho_u (z):=\olim\limits_{|s|\to +\infty\atop
s\in \IC}\ u(sz)-\log |s| \tag 1.16$$

Then $\rho(z)\in \cH$.

\head 2. Equilibrium measures\endhead
Let $E$ be a compact set in $\IC^N$ and $w$
an admissible weight function on $E$. We associate
the set $Z=Z(E, w)\subset \IC^{N+1}$ defined as follows:
$$Z:=\{(t,\lambda_1 t, \cdots,\lambda_N t)\in \IC^{N+1}\bigm|
(\lambda=\lambda_1, \cdots, \lambda_N)\in E,\ t\in \IC\ {\text{and}}\
|t|=w(\lambda)\} \tag 2.1$$
We will relate the weighted potential theory on $E$
with weight $w$ to potential theory on $Z$. We
will use the notation $(t, z)$ for a point in $\IC^{N+1}$
where $t\in \IC$ and $z\in \IC^N$. We denote by $C_\lambda$
the complex line in $\IC^{N+1}$ given by
$$C_\lambda:=\{(t, z)\in \IC^{N+1}| z_j=\lambda_j t\ {\text{for}}\
j=1, \cdots, N\ {\text{with}}\ (\lambda_1, \cdots, \lambda_N)\in \IC^N\
{\text{and}}\ \ t\in \IC\}\tag 2.2$$

$Z$ is a circular set, i.e., if $(t, z)\in Z$
then $(e^{i\theta} t, e^{i\theta} z)\in Z$ for all
$\theta \in [0, 2\pi)$. $\overline Z$ (the closure
of $Z$) is compact and circular.

$Z$ is non-pluripolar since $E$ is non-pluripolar
([BL2], lemma 6.1). Note that for $w\equiv 1$ (the
``unweighted" case) then $Z=\{|t|=1\}\times E$

Since $Z\subset \bigcup_{\lambda\in E} C_\lambda$ the same is true
for $\overline Z$ and so
$\overline Z\cap\{(t, z)\in \IC^{N+1} | t=0\}$ is either
empty (if $w$ is bounded below, on $E$, by a positive
constant) or else consists only of the origin.

\proclaim{Proposition 2.1}
$H_Z=H_{\bar Z}$ and $V_Z=V_{\bar Z}$
\endproclaim

\prf 
To prove the first statement it
suffices to show that if $u\in\cH(\IC\times \IC^N)$,
$u\le 0$ on $Z$ then $u\le 0$ on $\bar Z$. But since $w$
is u.s.c. we have, for $\lambda\in E$:
$$\overline Z\cap C_\lambda \subset\{(t, z)\in \IC_\lambda \bigm|
|t| \le w(\lambda)\}\tag 2.3$$
Applying the maximum principle to the
subharmonic function $t\to u (t,\lambda_1 t, \cdots, \lambda_N t)$ we
have $u\le 0$ on $\bar Z$.

The second statement follows similarly.

\proclaim{Proposition 2.2}
$V_Z=\Max (0, H_Z)$ and $\rho_{V_Z}=H_Z^*$.
\endproclaim

The first statement follows from ([Si1], proposition 5.6) and
the second from homogeneity (see [BLM], lemma 5.1).

\proclaim{Proposition 2.3}
$d\mu_{eq}(Z)$ has compact support in $\IC^{N+1}-\{t=0\}$.
\endproclaim

\prf
Since $H^*_Z\in \cH$ and $H^*_Z(0)=-\infty$
then $H^*_Z<0$ in a neighborhood of the
origin (estimates on the size of that neighborhood,
know as the Sibony-Wong inequality, can be found
in ([A], [Si2]).

Thus, the origin is an interior point of
${\hat\oZ}$. But $d\mu_{eq} (Z)$ places no mass on the
interior of ${\hat\oZ}$ so the result follows.\qed

Let $L$ denote the mapping
$L:\IC^{N+1}-\{t=0\}\to \IC^N$ given by
$$L(t, z)={z\over t}:=\lambda\in \IC^N\tag 2.4$$

If we consider $\IP^N$ (complex projective $N$-space)
as the space of lines through the origin in $\IC^{N+1}$,
then $L$ gives one of the standard coordinate charts for $\IP^N$.
Note that $L(C_\lambda)=\lambda$.

We recall the ``$H$-principle" of Siciak [Si3].
There is a natural $1-1$ correspondence between
$\cH(\IC^{N+1})$ and $\cL (\IC^N)$ as follows: To
$\tilde u(t, z)\in \cH (\IC\times \IC^N)$ associate
$$u(z):=\tilde u(1, z)\tag 2.5$$

Then $u\in \cL(\IC^N)$. Conversely, given 
$u\in \cL (\IC^N)$ we let
$$\tilde u (t, z):=u({z\over t})+\log|t| =L^* u(\lambda)+\log|t|\
{\text{for}}\ t\not= 0\tag 2.6$$
and
$$\tilde u(0, z)=\olim\limits_{|s|\to +\infty\atop
s\in\IC} u(sz)-\log (s)=\rho_u (z)\tag 2.7$$
Then $\tilde u\in \cH (\IC \times \IC^N)$.

Furthermore, let $P_d(t, z)$ be a homogeneous
polynomial of degree $d$ on $\IC\times \IC^N$. Then
${1\over d}\log |P_d (t, z)| \in \cH (\IC\times \IC^N)$ and
the associated
element of $\cL(\IC^N)$ under (2.5) is ${1\over d}\log|P_d(1, z)|$.
$P_d(1, z)$ is, of course a polynomial in $z$ of degree $\le d$.

Conversely, given a polynomical $G_d(z)$ in
$z$ of degree $\le d$, then ${1\over d} \log |G_d(z)|\in \cL(\IC^N)$.
The associated (via (2.6))
function in $\cH(\IC\times \IC^N)$ is ${1\over d}\log\bigm|
t^dG_d(z/t)\bigm|$.
Note that $t^d G_d({z\over t})$ is a homogeneous polynomial
on $\IC^{N+1}$ of degree $d$ in $(t, z)$. We use the notation:
$$P_d (t, z):=t^d G_d(z/t)\tag 2.8$$

Given a weighted polynomial $w^d G_d (\lambda)$ on $E$
we can relate its norm on $E$ with the norm of the
associated polynomial $P_d(t, z)$ on $Z$ (or equivalently, 
$\overline Z$).
Specifically, we have

\proclaim{Lemma 2.1}
$$\|w^d G_d \|_E =\|P_d (t, z)\|_Z$$
\endproclaim

\prf
For $(t, z)\in Z\cap C_\lambda$ and $\lambda \in E$, then
$$P_d(t, z)=t^d P_d(1, z/t)=t^d G_d (\lambda)\ {\text{so}}\tag 2.9$$
$$|P_d (t, z)|=|t|^d |P_d(1, z/t)|=w (\lambda)^d|G_d (\lambda)|
\tag 2.10$$
The result follows. \qed

Theorem 2.1 below gives the relation
between the pluricomplex Green function on $E$
and the homogeneous pluricomplex Green function of $Z$.

\proclaim{Theorem 2.1}
$H_Z(t, z)=L^* (V_{E, Q})+\log|t|$ for
$t\not=0$
\endproclaim

\prf
Let $\tilde u(t, z)\in \cH(\IC\times \IC^N)$ and suppose
$\tilde u\le 0$ on $Z$. Now
$\tilde u(t, \lambda_1 t, \cdots, \lambda_N t)=\tilde u(1, \lambda)
+\log|t|$ so for $\lambda\in E$, we have, 
$\log |t|+\tilde u(1, \lambda)\le 0$
on $C_\lambda \cap Z$.
Thus,
$$u(\lambda)\le -\log |t|=Q(\lambda)\ {\text{on}}\ C_\lambda\cap Z
\tag 2.11$$

Hence $u\le V_{E,Q}$ and so, using (2.6)
$$\tilde u(t, z)\le L^*(V_{E,Q})+\log |t| \ {\text{for}}\ t\not=0$$.

Taking the pointwise sup in $(t, z)$ over all
such $\tilde u$ we have
$$H_Z\le L^* (V_{E,Q})+\log |t| \ {\text{for}}\ t\not=0\tag 2.12$$

It remains to prove the reverse inequality. Given 
$u \in \cL(\IC^N)$ with $u\le Q$ on $E$ then for
$(t, z)\in C_\lambda\cap Z$, using (2.6), we have
$\tilde u(t,z)=\log |t|+u(\lambda)\le \log |t|-\log w(\lambda)\le 0$.
Hence $\tilde u (t, z)\le 0$ on $Z$ and
$\tilde u(t, z)\le H_Z$ on $\IC\times\IC^N$. That is,
$$L^*(u)+\log |t| \le H_Z\ \ {\text{for}}\ \ t\not=0\tag 2.13$$

Taking the pointwise sup over all such $u$ gives
the reverse inequality to (2.12).\qed

\proclaim{Corollary 2.1}
$H^*_Z (t, z)=\widetilde{V^*_{E, Q} (t, z)}$
\endproclaim

\prf
Consider a point $(t_0, z_0)$ with $t_0\not=0$. 
Let $\lambda_0:={z_0\over t_0}$

Then, by theorem 2.1,
$$\olim\limits_{(t, z)\to (t_0, z_0)} H_Z (t, z)=
\olim\limits_{\lambda \to \lambda_0} V_{E, Q} (\lambda)+
\log |\lambda_0|$$
and the right side is $\widetilde{V^*_{E, Q} (t_0, z_0)}$ by (2.6).
This proves the result for $t\not=0$ but since both
sides (in the statement of corollary 2.1) are PSH
functions on $\IC^{N+1}$ and agree for $t\not=0$ they
must agree on $\IC^{N+1}$. \qed

Note that the result in ([K], prop. 2.9.16) is
similar but not immediately applicable.

Corollaries 2.2, 2.3 and 2.4 deal with the converge of
sequences of pluricomplex Green functions for sequences
of weights converging in various manners (see also
lemma 7.3 [BL2]).

\proclaim{Corollary 2.2}
Let $E\subset \IC^N$ be
compact and $\{w_j\}_{j=1,2, \cdots}$ a sequence of admissible
weights on $E$ and let $w$ also be an admissible weight on $E$.
Suppose that $w_j \downarrow w$. Then
$$\lim\limits_j V_{E, Q_j}=V_{E, Q}.$$
\endproclaim

\prf
Let $Z_j:=Z_j (E, w_j)$ and $Z:=Z(E, w)$  be
the associated circular sets in $\IC^{N+1}$. But 
$\{(t, z)\bigm| |t|\le w_j\}\downarrow \{(t, z)\bigm| |t|\le w\}$. Hence
$\hat{\bar Z}_j \downarrow \hat{\bar Z}$ and the result follows from
theorem 2.1 and ([K], corollary 5.1.2).

\proclaim{Corollary 2.3}
Let $E, \{w_j\}, w$ be as in 
corollary 2.2, except that
$w_j\downarrow w$ q.e. then
$$\lim\limits_j V^*_{E, Q_j}=V^*_{E, Q}$$
\endproclaim

\prf
$\bigcap_j\{(t, z)\bigm| |t|\le w_j\}$ and
$\{(t, z)\bigm| |t|\le w\}$ differ
by a pluripolar set so the result follows from
([K], cor. 5.2.5).

\proclaim{Corollary 2.4}
Let $E, \{w_j\}$, $w$ be as in corollary 2.2
except that $w_j \uparrow w$ q.e. Then 
$\lim\limits_j V^*_{E, Q_j}=V^*_{E, Q}$.
\endproclaim

\prf
For some pluripolar set $F$ we have
$Z_j \cup F\uparrow Z\cup F$. Hence by ([K], cor. 5.2.5 and 5.2.6)
$$V^*_{Z_j} =V^*_{Z_j \cup F}\downarrow V^*_{Z\cup F}=V^*_Z$$
Hence, using homogeneity, $H^*_Z \downarrow H_Z$ and by
corollary 2.1,
$$V^*_{E, Q_j}\downarrow V_{E, Q}.\qquad\qquad\qquad\qquad\qed$$

By proposition 2.3, we may consider
$L_* (d\mu_{eq} (Z))$-the push forward of the
equilibrium measure $d\mu_{eq} (Z)$ under $L$. Since
$\supp (d\mu_{eq} (Z))\subset Z\subset\bigcup_{\lambda\in E}
\IC_\lambda$ we have $\supp(L_* (d\mu_{eq}(Z))\subset E$.
There is however a more precise relation. Assume that
$Z$ is regular. 
The equilibrium measure on $Z$ and the weighted equilibrium
measure on $E$ are related by:

\proclaim{Theorem 2.2}
$L_* \Big({1\over 2\pi} d\mu_{eq} (Z)\Big)=d\mu_{eq} (E, w)$.
\endproclaim

\prf
The proof is based on lemma 3.3 in [DeM] which
itself is based on work of Briend. (Note that
we use the convention of Klimek's book [K] for
$d^c:=i(\bar\partial -\partial)$ not that of [DeM]. This results in the
factor ${1\over 2\pi}$ in the statement of theorem 2.2).

$H_Z$ is continuous by proposition 2.2 so, as a consequence of Theorem 2.1
$V_{E, Q}$ is continuous. Then
$$dd^c H_Z=dd^c L^*(V_{E, Q})=L^*(dd^c V_{E, Q})\ \
{\text{for}}\ \ t\not=0\ \ {\text{and so}}\tag 2.14$$
$$(dd^cH_Z)^N=L^*(dd^c V_{E, Q})^N=L^*(d\mu_{eq}(E, w))\
{\text{for}}\ \  t\not= 0.\tag 2.15$$

Let $\phi$ be a smooth compactly
supported function on $\IC^{N+1}-\{t=0\}$. then
$$\int\limits_{\IC^{N+1}}\!\! \phi d\mu_{eq} (Z)=
\int\limits_{\IC^{N+1}}\!\! \phi (dd^c V_Z)^{N+1} =
\int\limits_{\IC^{N+1}}\!\! V_Z dd^c \phi \wedge (dd^c V_Z)^N
\tag 2.16$$

Now, as $\ep\downarrow 0$,
$\Max (H_Z, \ep)-\ep \uparrow V_Z$ uniformly on $C^N$.
so
$$\int\limits_{\IC^{N+1}}\!\! V_Z dd^c\phi\wedge (dd^c V_Z)^N=
\lim\limits_{\ep \to 0}\int_{\IC^{N+1}}\!\!
(\Max (H_Z, \ep)-\ep) dd^c\phi\wedge
(dd^c V_Z)^N\tag 2.17$$
Note that $\{z\in\IC^{N+1}\bigm| \Max(H_Z, \ep)=\ep\}$ is
a neighborhood of $\bar Z$ in $\IC^{N+1}$ and if
$\Max (H_Z (z), \ep)>\ep$ then $V_Z=H_Z$. So, in the
expression on the right of (2.17) we may replace
$V_Z$ by $H_Z$ to obtain.

$$\int\limits_{\IC^{N+1}}\!\! V_Z dd^c\phi\wedge (dd^c H_Z)^N =
\int\limits_{\IC^{N+1}}\!\! V_Z dd^c \phi\wedge L^* (dd^c V_{E,Q})^N
\tag 2.18$$

The integrands in (2.16), (2.17) and (2.18) all have
compact support in $\IC^{N+1}-\{t=0\}$ so
the right side of (2.18) is equal to
$$\int\limits_{\IC^N}\!\! \Big(\int\limits_{C_\lambda}
V_Z dd^c \phi\Big) d\mu_{eq} (E, w)=
\int\limits_{\IC^N}\!\! \Big(\int\limits_{\IC_\lambda}
\phi dd^c V_Z\Big) d\mu_{eq} (E, w)\tag 2.19$$
For $\lambda\in E$, we let $dm_\lambda$ be the Lebesgue measure on the
circle $|t|=w(\lambda)$ in $C_\lambda$ normalized to have total mass 1.
Then ${1\over 2\pi} dd^c V_{Z/{C_\lambda}}=dm_\lambda$.

The right side of (2.19) is thus equal to
$${1\over 2\pi} \int\limits_{\IC^N}\!\!
\Big(\int\limits_{C_\lambda}\!\! \phi dm_\lambda\Big) d\mu_{eq}
(E, w)\tag 2.20$$
which proves theorem 2.2.\qed

The next corollary shows that
the assumption of $\bar Z$ being regular may be dropped
from the hypothesis of theorem 2.2.

\proclaim{Corollary 2.5}
$L_*\Big({1\over 2\pi} d\mu_{eq} (Z)\Big)=d\mu_{eq} (E, w)$.
\endproclaim

\prf
We need only find a sequence of
locally regular compact sets $E_j$, admissible, continuous, weights
$w_j$ on $E_j$ such that $E_j\downarrow E$ and
$w_j \downarrow w$ on $E$.
Then $Z (E_j, w_j)\downarrow Z(E, w)$. Applying theorem 2.2
to each $Z(E_j, w_j)$ and $E_j$ and taking limits gives
the result.

To construct such a sequence of $E_j$ and $w_j$
we may follow the procedure of ([BL2], section 7).

\head 3. The Bernstein-Markov inequality\endhead
Given a compact set $E\subset \IC^N$ and a finite
positive Borel measure $\mu$ on $E$, we say that
$(E, \mu)$ satisfies the Bernstein-Markov (B-M) inequality
if, for every $\ep >0$, there exists a constant $C=C(\ep)>0$
such that, for all holomorphic polynomials $p$ we have
$$\|p\|_E \le C (1+\ep)^{\deg (p)} \|p\|_{L^2 (\mu)}\tag 3.1$$
This inequality may be used to relate $L^2$
properties of polynomials with potential
theoretic invariants of $E$ (see [B1] and [BL2] for
conditions under which the inequality holds).

We will introduce a ``weighted" 
version of the B-M inequality.

Given a compact set $E\subset \IC^N$, an admissible
weight $w$ on $E$ and a finite positive Borel
measure $\mu$ on $E$, we say that $(E, w, \mu)$
satisfies the weighted B-M inequality if for
all $\ep>0$, that exists a constant $C=C(\ep)>0$ such that,
for all weighted polynomials $w^d p$ we have
$$\|w^d p\|_E \le C(1+\ep)^d \| w^d p\|_{L^2(\mu)}\tag 3.2$$
Of course, for $w\equiv 1$, (3.2) reduces to (3.1).

We will relate the weighted B-M inequality for 
$(E, w, \mu)$ to a B-M inequality
on $\bar Z$ with respect to a certain  associated measure $\nu$.
The measure $\nu$ is defined as follows:
$$d\nu=dm_\lambda\otimes d\mu\ {\text{for}}\ \
\lambda \in E\ \ {\text{so\ that}}\ \ \supp(\nu)\subset
\bigcup_{\lambda\in E} C_\lambda\tag 3.3$$
That is, for $\phi$ continuous with compact support
in $\IC^{N+1}-\{t=0\}$ we have
$$\int\limits_{\IC^{N+1}}\!\! \phi d\nu=
\int\limits_E\!\!\Big(\int\limits_{C_\lambda}\!\!
\phi dm_\lambda\Big) d\mu(\lambda)\tag 3.4$$

\proclaim{Theorem 3.1}
$(E, w, \mu)$ satisfies the weighted B-M inequality
if and only if $(\bar Z, \nu)$ satisfies the B-M inequality.
\endproclaim

\prf
First (using the notation of lemma 2.1) we prove

\proclaim{Lemma 3.1}
$$\|P_d\|_{L^2(\nu)}=\|w^d G_d\|_{L^2 (\mu)}$$
\endproclaim

\prf
$P_d$ is a homogeneous polynomial of
degree $d$ on $\IC^{N+1}$.
Now
$$\align
\int\limits_{\IC^{N+1}}\!\! |P_d (t, z)|^2 d\nu & =
\int\limits_E\!\! \Big(\int\limits_{C_\lambda}\!\!
|P_d (t, z)|^2 dm_\lambda\Big) d\mu (\lambda) \\
&=\int\limits_E\!\! |w^d G_d|^2 d\mu (\lambda),\ {\text{using}}\
(2.9) \\
&=\|w^d G_d\|_{L^2(\mu)}\qquad\qquad\qquad \qed\endalign$$

Now, suppose $(\bar Z, \nu)$ satisfies the B-M inequality.
Applying that inequality to homogeneous polynomials,
using lemmas 2.1 and 3.1 we obtain the
weighted B-M inequality for $(E, w, \mu)$.

For the converse, suppose $(E, w, \mu)$ satisfies
the weighted B-M inequality.

We first note that if two monomials are
of different degrees, they are orthogonal in $L^2 (\nu)$
since their restrictions to any $C_\lambda$ are orthogonal in
$L^2 (dm_\lambda)$. Hence for a polynomial $p$ on $\IC^{N+1}$,
written as a sum of homogeneous polynomials
$$p=\sum^d_{i=0} p_i\ \ {\text{then}}\tag 3.5$$
$$\|p\|_{L^2(\nu)} =\sum^d_{i=0} \|p_i\|_{L^2(\nu)}\tag 3.6$$

Hence, for any $\ep >0$ there is a $C>0$ such that
$$\|p\|_Z \le \sum^d_{i=0} \|p_i\|_Z \le C(1+\ep)^d
\|p_i\|_{L^2(\nu)} \le C(d+1) (1+\ep)^d\|p\|_{L^2(\nu)}\tag 3.7$$
where the second inequality comes from lemmas 2.1, 3.1
and the weighted B-M inequality for $(E, w, \mu)$. The
third inequality in (3.7) comes from (3.6). The B-M
inequality for $(\bar Z, \nu)$ follows from (3.7).\qquad \qquad \qed

\proclaim{Coroallary 3.1}
Suppose $\bar Z$ is regular. Then $(E, w, d\mu_{eq} (E, w))$
satisfies the weighted B-M inequality.
\endproclaim

\prf
It is a result of Nguyen-Zeriahi [NZ]
combined with ([K], corollary 5.6.7) that $(\bar Z, d\mu_{eq} 
(\bar Z))$
satisfies the B-M inequality. However, by
theorem 2.2,
${1\over 2\pi} d\mu_{eq} (\bar Z)=dm_\lambda \otimes
d\mu_{eq} (E, w)$.\qed

We will give another general situation in which
the weighted B-M inequality holds (see also [StTo], theorem
3.2.3. (vi)).

\proclaim{Theorem 3.2}
Let $E$ be a locally regular,
compact set $\subset \IR^N$ and let $w$ be continuous on
$E$ with $\inf_{z \ep E} w(z)>0$.
Suppose that $\sigma$ is a finite positive Borel measure on $E$
and $(E, \sigma)$
satisfies the B-M inequality. Then $(E, w, \sigma)$ satisfies the
weighted B-M inequality.
\endproclaim

\prf
$\log w$ is continuous on $E$ and so may be,
by the Weierstrass theorem, approximated by (real)
polynomials. That is, given $\ep>0$ there exists
$g_\ep=g_\ep (x_1,\cdot, x_n)$, a real polynomial, such that
$\|\log w-g_\ep\|_E \le \ep$
Taking exponentials, we have
$$e^{-\ep}\le {w\over \exp (g_\ep)} \le e^\ep\quad {\text{for}}\
z\in E\tag 3.9$$

We consider $g_\ep$ as a holomorphic polynomial
$$g_\ep =g_\ep (z_1, \cdots, z_N)$$

Taking sufficent many terms in the power 
series for $\exp (g_\ep)$ we get a holomorphic polynomial $H$
such that, for $\ep$ sufficiently small
$$ 1-2\ep \le {w\over H}\le 1+2\ep\quad {\text{for}}\ \
z\in E\tag 3.10$$

Now, consider a weighted polynomial $w^d G$ and let
$$J:=GH^d\tag 3.11$$
Then
$$w^d |G|=|J| ({w\over |H|})^d\quad {\text{so that}}\tag 3.12$$
$$|J|(1-2\ep)^d \le w^d|G| \le |J|(1+2\ep)^d\ \ {\text{for}}\ \
z\in E\tag 3.13$$

Now, by the B-M inequality for $(E, \sigma)$ we have,
given $\ep_1 >0$ a constant $C_1 > 0$ such that
$$\|J\|_E \le C_1 (1+\ep_1)^{(h+1)d} \|J\|_{L^2 (\sigma)}\tag 3.14$$

where $h:=\deg H$.
Hence
$$\align
\|w^dG\|_E & \le \|J\|_E (1+2\ep)^d\quad {\text{by\ the\ right\
inequality\ in\ (3.13)}} \\
&\le C_1 (1+\ep_1)^{(h+1)d} (1+2\ep)^d \|J\|_{L^2(\sigma)}\quad 
{\text{by\ (3.14)}} \\
&\le C_1(1+\ep_1)^{(h+1)d} {(1+2\ep)^d\over (1-2\ep)^d}
\|w^d G\|_{L^2(\sigma)}\quad  {\text{by}}\endalign$$
the left inequality in (3.13).

Now, $\ep>0$ having been chosen, and $h$ fixed,
we choose $\ep_1$ so that $(1+\ep_1)^{(h+1)} \le 1+\ep$ and so
we obtain
$$\|w^d G\|_E \le C_1{(1+\ep)^d (1+2\ep)^d\over(1-2\ep)^d}
\|w^d G\|_{L^2(\sigma)}\tag 3.16$$

But $\ep>0$ is arbitrary so the weighted B-M
inequality holds.

\noindent
{\bf Example 3.1}\
Let $B_R=\{x \in \IR^N\bigm| |x|\le R\}$ be the (real)
ball of radius $R$ (center the origin). Then $(B_R, dx)$
satisfies the BM inequality (see [B]) where $dx$ denotes
Lebesgue measure.

Let $w(x)$ be any continuous positive function on
$B_R$. Then by theorem 3.2 $(B_R, w, dx)$ satisfies the
weighted B-M inequality.

\head 4. $L^2$ theory of weighted polynomials\endhead
Let $E$ be a compact non-pluripolar
subset of $\IC^N$, $w$ an admissible
weight on $E$, and $\mu$ a finite positive Borel measure
with $\supp(\mu)=E$. For $d$ a positive integer,
the monomials are linearly independent in $L^2(w^{2d} \mu)$
([Bl1], prop. 3.5 adapts to this situation). Ordering
via a lexicographic ordering on
their multi-index exponents and applying the
Gram-Schmidt procedure we obtain orthonomal
polynomials $\{p_\alpha^d (z, \mu)\}_{\alpha\in  \IN^n}$.
They satisfy
$$\int\limits_{\IC^N}\!\! p^d_\alpha (z,\mu)p^d_\beta
(z, \mu)w^{2d} d\mu=\delta_{\alpha \beta}\tag 4.1$$
for $\alpha, \beta$ multi-indices.


We can write
$$p^d_\alpha (z,\mu) =a^d_\alpha z^\alpha +\ \
{\text{(monomials\ of\ lower\ lexicographic\ order)\ where}}\
a^d_\alpha>0 \tag 4.2$$
\vskip -.2truein
\noindent
We will only consider these polynomials where 
$|\alpha|=d$.

In the case that $(E, w, \mu)$ satisfies the weighted B-M
inequality we will show that the leading exponents
$\{a^d_\alpha\}$ have asymptotic limits in the following sense.
First we let
$$\Sigma_0:=\{\theta \in\IR^N \bigm| \theta=(\theta_1, \cdots,
\theta_N),\ \sum^n_{j=1} \theta_j=1,\ \theta_j>0\}\tag 4.3$$

We consider sequences of multi-indices $\{\alpha (j)\}$
with, for some $\theta\in \Sigma_0$
$$\lim\limits_j |\alpha(j)|=+\infty\ \ {\text{and}}\ \
\lim\limits_j {\alpha(j)\over |\alpha (j)|}=\theta\tag 4.4$$

\proclaim{Theorem 4.1}
Suppose $(E, w, \mu)$ satisfies the weighted B-M
inequality and $\{\alpha(j)\}$ is a sequence of multi-indices
satisfying (4.4). Then
$$\lim\limits_j 
\Big(a_{\alpha(j)}^{d}\Big)^{1\over d}=
{1\over \tau^w (E, \theta)}$$
where $\tau^w (E, \theta)$ is the weighted directional Tchebyshev
constant of $E$ in the direction $\theta$.
\endproclaim

\prf
First, we recall the definition of weighted
directional Tchebyshev constant (see [BL2]). For $\alpha$ a multiindex
we let $P(\alpha)=\{q|q =
z^\alpha+\sum_{\beta<\alpha} c_\beta z^\beta\}$
where $c_\beta\in \IC$ and the notation $\beta<\alpha$ is used to denote
the fact that the multiindex $\beta$ preceeds $\alpha$ in the
lexicographic ordering on the multi-indices.

For $\alpha$ a multiindex with $|\alpha|=d$ we let
$t^d_\alpha$ denote a (Tchebyshev) polynomial which
minimizes $\{\|w^d q\|_E \bigm| q\in \cP(\alpha)\}$. That is,
$t^d_\alpha \in \cP (\alpha)$ and 
$$\|w^d t^d_\alpha\|_E = \inf \{\|w^d q\|_E \bigm| q\in\cP (\alpha)\}.
\tag 4.5$$

Then (see [BL2]) it is known that for a 
sequence of multi-indices $\{\alpha(j)\}_{j=1, 2, \cdots}$
satisfying (4.4) the limit
$$\tau^w (E, \theta):=\lim\limits_j
\|w^d t^d_{\alpha(j)}\|^{1\over d}_E\tag 4.6$$
exist and is called the weighted 
Tchebyshev constant in the direction $\theta\in\Sigma_0$.

Now, it follows from general Hilbert space theory that
$$a^d_\alpha ={1\over \|w^d q^d_\alpha\|_{L^2(\mu)}}\tag 4.7$$
where $q^d_\alpha$ is the unique polynomial in $\cP(\alpha)$
satisfying.
$$\|w^d q^d_\alpha\|_{L^2(\mu)}=\inf\{\|w^d q\|_{L^2(\mu)} 
\bigm| q\in\cP(\alpha)\}\tag 4.8$$

Now, for $\ep>0$, there is a $C>0$ such that
$$\aligned
\|w^d q^d_\alpha \|_E & \le
C(1+\ep)^d \|w^d q^d_\alpha\|_{L^2(\mu)}\quad 
{\text{by\ the\ weighted\ B-M\ inequality}}\\
& \le C(1+\ep)^d \|w^d t^d_\alpha\|_{L^2(\mu)}\quad
{\text{by\ (4.8)}}\\
& \le C_1(1+\ep)^d \|w^d t^d_\alpha\|_E \endaligned\tag 4.9$$
since the sup norm estimates the $L^2$ norm for a finite
measure with compact support.

Hence, for every $\ep>0$ there is a constant $C_1>0$ such that
$$\|w^d t^d_\alpha\|_E \le \|w^d q^d_\alpha\|_E\le
C_1(1+\ep)^d \|w^d t^d_\alpha\|_E
\tag 4.10$$

Now, given a sequence of multi-indices $\{\alpha(j)\}$
satisfying (4.4), taking the $1/d$ powers of the 
expressions in (4.10), letting $j\to \infty$, using
(4.6), (4.7) and the fact that
$\ep > 0$ is arbitrary, the result follows.\qed

{\bf Example 4.1}\
On $\IR^N$ we consider orthonormal
polynomials with respect to the inner product
given by $e^{-H(x)} dx$ where $dx$ is Lebesque measure on $\IR^N$
and $H(x)$ satisfies
$$\aligned
i)\quad & H(x)\ {\text{is\ homogeneous\ of\ degree}}\ \gamma>0.\
{\text{That\ is}}\\
& H(cx)=c^\gamma H(x)\quad c\in\IR.\\
ii)\quad & H(x)>0\quad {\text{for\ all}}\ \ x\not=0\endaligned
\tag 4.11$$

We let $\{p_\alpha (x)\}_{\alpha\in \IN^N}$
denote the orthonormal
polynomials obtained, by applying the
Gram-Schmidt procedure to the (real) monomials
ordered via a lexicographic ordering of their exponents. Then
$$\int\limits_{\IR^N}\!\! p_\alpha (x) p_\beta (x) e^{-H(x)}
dx=\delta_{\alpha\beta}\tag 4.12$$
for any two multi-indices $\alpha, \beta$.

We write
$$p_{\alpha}(x)=a_\alpha x^\alpha +\ {\text{(sum\ of\ monomials\ of\
lower\ lexicographic\ order).}}\ a_\alpha>0 \tag 4.13$$ 

We will obtain asymptotic estimates (see (4.24))
for $|a_\alpha|^{1\over |\alpha|}$
for a sequence of multi-indices satisfying (4.4).
In the case $N=1$, these estimates are Theorem VII, 1.2 
of [SaTo]. In that case explicit knowledge of the set
$S_w$ (defined below) yields an explicit form to the
right hand side of (4.24). It would be of interest to 
find $S_w$ explicitly in the case $N>1$.

In the one-dimensional case $(N=1)$ this
gives a version of so-called weak asymptotics and
in this case considerably more detailed asymptotic
results are known (see [SaTo] or [Dei]).

>From general Hilbert space thoery,
$$a_\alpha^{-1}=\inf \{\|e^{-{H(x)\over 2}} q(x)\|_{L^2(\IR^N)}
\bigm| q\in \cP_\IR (\alpha)\}\tag 4.14$$
where $q\in \cP_\IR (\alpha)=$ \{polynomials of the form
$x^\alpha+\sum_{\beta<\alpha} r_\beta x^\beta$ with 
$r_\beta\in \IR$\}.

For $| \alpha|=d$ we scale by $x=d^{1\over\gamma} y$. We get
$$a^{-1}_{\alpha} =d^{d\over \gamma}\ \inf 
\Big\{\|e^-{dH(y)\over 2}
q(y)\|_{L^2(\IR^N)} \bigm| q\in \cP_R(\alpha)\Big\}.\tag 4.15$$
Consider the weight $w(y)=e^{-{H(y)\over 2}}$
on $\IR^N\subset \IC^N$. This weight is admissible in the
sense of ([SaTo], appendix B) although, since $\IR^N$ is not
compact, not in the sence of 1.10. We let $Q(y)={H(y)\over 2}$.
The following is known ([SaTo], appendix B).
$S_w:=(dd^c V_{\IR^N, Q})^N$ has compact support. For any
weighted polynomial $w^dp$ we have
$$|w^d p(y)|\le \|w^d p\|_{S_w} \exp (d(V_{\IR^N, Q} - Q))\tag 4.16$$
In particular
$$\sup_{\IR^N} |w^d p|=\sup_{S_w} |w^d p|.\tag 4.17$$

Now $V_{\IR^N, Q}\in \cL$ so, fixing $R>0$ large, using (4.16), 
there is
a constant $A>0$ such that.
$$|w^d p(y)|\le \|w^d p\|_{S_w} e^{-Ad|y|^\gamma}\ {\text{for}}\
|y|\ge R\tag 4.18$$

Now
$$\|w^dp\|^2_{L^2(\IR^N)} \le \|w^dp\|^2_{L^2(B_R)}+
\|w^d p\|^2_{S_w} \int\limits_{|y|\ge R}\!\!
e^{-2Ad|y|^\gamma} dy\tag 4.19$$

We may assume $S_w \subset B_R$. Then
$$\aligned
&\|w^dp\|_{S_w}  =\|w^d p\|_{B_R} \le C (1+\ep)^d\|w^d p\|_{L^2(B_R)}\ \
{\text{since}} \\
& (B_R, w, dx)\ \ {\text{satisfies\ the\ weighted\ B-M\ inequality\
(see\ example\ 3.1)}}\endaligned \tag 4.20$$
Now simple estimates show there is a constant $c_1>0$
such that
$$\int\limits_{|y|\ge R}\!\! e^{-2Ad|y|^\gamma} dy
\le e^{-dc_1}\tag 4.21$$
We get
$$\|w^d p\|_{L^2(\IR^N)} \le \|w^d p\|_{L^2(B_R)}
\Big(1+{C^2(1+\ep)^{2d}\over e^{dc_1}}\Big)^{1\over 2}.\tag 4.22$$
However for $\ep>0$ sufficiently small, the expression on
the right of (4.22) is bounded in $d$.
It follows that for a sequence of multi-indices $\{\alpha(j)\}$
satisfying (4.4)
$$\lim\limits_{j\to\infty}[ \inf\{\|w^d q\|_{L^2(\IR^N)}\bigm|
q\in \cP_\IR (\alpha)\}^{1\over d}] =\lim\limits_{j\to\infty}
[\inf \{\|w^d q\|_{L^2(B_R)}\bigm| q\in\cP_\IR 
(\alpha(j))\}^{1\over d}]\tag 4.23$$
But, by (the proof of) theorem 4.1 the limit
on the right side of (4.23) exists and it may be
identified with $\tau^w(S_w, \theta)$ using (4.1).
Hence we obtain
$$\lim\limits_{j\to \infty} a_{\alpha(j)}^{1\over d} 
d^{1\over \gamma}={1\over \tau^w(S_w, \theta)}\tag 4.24$$

\head 5. General Circular Sets\endhead
The circular sets which arise in the form
$\hat{\bar Z} (E, w)$ (i.e. the polynomially convex hull
of a set of the form $Z(E, w)$) are
\item{i)\ } polynomially convex
\item{ii)\ } circular
\item{iii)\ } compact
\item{iv)\ } non-pluripolar

However, they are not the most general
sets with the properties i), ii), iii), iv). We will show,
however,
that the most general set with those properties is,
in an appropriate sense, a limit of sets of the
form $\hat{\bar Z} (E, w)$

Let $Z\subset \IC^{N+1}$ be a set with properties 
i), ii), iii), iv) above.
Then the origin is an interior point of $Z$. We associate
to $Z$ the a function on $\IC^N$ defined by
$$w(\lambda):=\sup \{|t| \bigm| (t, z) \in Z \cap C_{\lambda}\}
\tag 5.1$$

Then $w(\lambda)>0$ for  all $\lambda$ and $w$ is bounded above. We let
$Q(\lambda):=-\log w(\lambda)$.

\proclaim{Proposition 5.1}
$w$ is u.s.c. on $\IC^N$.
\endproclaim

\prf
Fix $\lambda^0 \in \IC^N$. Let $\{\lambda^s\}_{s=1, 2, \cdots}$
be a 
sequence in $\IC^N$ converging to $\lambda^0$. We may suppose,
passing to a subsequence if necessary, that
$\lim\limits_s (\lambda^s):=w^0$ exists.
The points $(w (\lambda^s),\ \lambda^s_1 w(\lambda^s), \cdots,
\lambda^s_N w(\lambda^s))\in Z\cap C_{\lambda^s}$ so,
since $Z$ is compact, the point $(w^0, \lambda^0_1 w^0, \cdots,
\lambda^0_N w^0)\in Z\cap C_{\lambda^0}$.
Thus, by definition of $w$, $w(\lambda^0)\ge w^0=\lim\limits_s
w (\lambda^s)$. \qquad\qquad \qed

\proclaim{Proposition 5.2}
$Q^*$ is PSH on $\IC^N$.
\endproclaim

\prf
$Z=\{(t, z) \in \IC^{N+1} \bigm| H_Z (t, z)\le 0\}$

Now
$$\aligned
H_Z(t, z) &=\log|t|   + H_Z (1, \lambda)\ \ {\text{so that}}\\
\log w(\lambda) &=-H_Z (1, \lambda)\ \ {\text{and}}\\
Q(\lambda)&= H_Z (1, \lambda)\endaligned$$

But $H_Z=H_Z^*$ outside a circular pluripolar set in $\IC^{N+1}$
so by [BL2, lemma 6.1] $H_Z(1, \lambda) =H^*_Z (1, \lambda)$
q.e. on $\IC^N$.

Thus $Q^*=H^*_Z (1, \lambda)$ \qquad\qquad \qed.

\noindent
{\bf Example 5.1}\
$Z=\{|t|^2 + |z_1|^2 +\cdots +|z_N|^2 \le 1\}\subset \IC^{N+1}$.
Then $w(\lambda)=(1+|\lambda_1|^2+\cdots + |\lambda_N|^2)^{-{1\over2}}$
and
$Q(\lambda)={1\over 2} \log (1+ |\lambda_1|^2 +\cdots
| \lambda_N|^2)$.

As the above example illustrates, in general the
functions $w(\lambda)$ which arise in this way
are {\it{not}} admissible weights in the sense of
[SaTo], appendix B, definition 2.1. In particular, 
$\lim\limits_{|\lambda| \to \infty} |\lambda|w(\lambda)\not=0$

Let $Z_R:=\{(t, z)\in Z\bigm| {|z| \over |t|}\le R\}$.

\proclaim{Propostion 5.3}
$\lim\limits_{R\to\infty} d\mu_{eq} (Z_R)=d\mu_{eq} (Z)$
weak$^*$
\endproclaim

\prf
$Z_R \cup T\uparrow Z\cup T$ where $T$ is the pluripolar
set $\{(t, z)\in \IC^{N+1} \bigm| t=0\}$. Hence
$V^*_{Z_R} \downarrow V^*_{Z}$ as $R\to \infty$.
The result follows from
([K], cor 5.2.5 and 5.2.6) and the continuity of the
Monge-Amp\`ere operator under decreasing limits [K].\qquad\qquad
\qed

\proclaim{Proposition 5.4}
$\lim\limits_{R\to \infty} L_* ({1\over 2\pi}d\mu_{eq} (Z_R))=
{(dd^c Q^*)^N}$
\endproclaim

\prf
Let $w_R=w/_{B(0, R)}$
Let $Q_R=-\log w_R$. Then the family
of PSH functions $V^*_{B(0, R), Q_R}$ 
is decreasing as $R\to \infty$.
Since $Q^*$ is PSH, $V^*_{B(0, R), Q_R}
=Q^*$ on
$B(0, R)$ so $V^*_{B(0, R), Q_R} \downarrow Q^*$.
The result follows using theorem 2.2.\qquad\qquad \qed.
\bigskip
\vfill\eject

\vglue .5truein

Thomas Bloom

Department of Mathematics

University of Toronto

Toronto, ON

CANADA\ \ M5S 2E4
\bigskip

email: bloom\@math.toronto.edu

\end{document}